\documentclass{article}
\usepackage{amsfonts,amsmath,amssymb,latexsym}
\newcounter{theorem}

\newcounter{lemma}

\newcounter{proposition}

\newcounter{example}

\newcounter{corolary}

\begin{document}

\begin{center}
{\bf \Large The continuous extension of the logarithmic double layer potential to the Ahlfors-regular boundary}

\vskip 2mm

{\bf \large Sergiy Plaksa}

\vskip 2mm

{Institute of Mathematics of NAS of Ukraine, Kyiv, Ukraine\\
University of Padova, Padova, Italy}

\vskip 2mm

{plaksa62@gmail.com}
\end{center}

\begin{abstract}
\noindent For the real part of the Cauchy-type integral that is known to be the logarithmic potential of the double layer, 
a necessary and sufficient condition for the continuous extension to the Ahlfors-regular boundary is established.\\[1mm]
\end{abstract}

\section{Introduction}
\label{Introduction}

Let $\gamma$ be a closed rectifiable Jordan curve in the complex plane $\mathbb{C}$, and
let $D^{+}$ and $D^{-}$ be the interior and exterior domains bounded by $\gamma$, respectively.

The classical theory of the logarithmic double layer potential
\begin{equation}\label{log-pot}
\frac{1}{2\pi}\int\limits_{\gamma} g(t)\frac{\partial}{\partial {\bf n}_t}\left(\ln\frac{1}{|t-z|}\right)\,ds_t
 \qquad \forall\,z\in D^{\pm}
\end{equation}
is developed in the case where the integration curve $\gamma$ is a Lyapunov curve (see, for example, J.~Plemelj \cite{Plemelj-1911}).
Here ${\bf n}_t$ and $s_t$ denote the unit vector of the outward normal to the curve $\gamma$ at a point $t\in\gamma$ and an arc coordinate of this point, respectively, and the integral density $g : \gamma\rightarrow\mathbb{R}$ takes values in the set of real numbers $\mathbb{R}$.

J.~Radon \cite{Radon-46} established the continuous extension of the logarithmic double layer potential from the domains $D^{+}$ and $D^{-}$ to the boundary $\gamma$ in the case where $\gamma$ is a curve of bounded rotation, i.e., a curve for which the angle between the tangent to the curve and a fixed direction is a function of bounded variation.
It is known that the class of Lyapunov curves and the class of Radon curves of bounded rotation are different, i.e., each of them contains curves that do not belong to the other class (see, for example, I.I.~Danilyuk \cite[p.~26]{Dan-m}).

J.~Kr\'al \cite{Kral-1-1964} established a necessary and sufficient condition for the curve $\gamma$,
under which the logarithmic double layer potential is continuously extended from the domains
$D^{+}$ and $D^{-}$ to the boundary $\gamma$ for all continuous functions $g : \gamma\rightarrow\mathbb{R}$.

The logarithmic double layer potential (\ref{log-pot}) is the real part of the Cauchy-type integral
(see, for example, F.D.~Gakhov \cite{Gah}, N.I.~Muskhelishvili \cite{Mus})
\begin{equation}\label{C-type-int}
\widetilde{g}(z):=
\frac{1}{2\pi i}\int\limits_{\gamma}\frac{g(t)}{t-z}\,dt
 \qquad \forall\,z\in D^{\pm}\,.
\end{equation}

The theory of the boundary properties of the integral (\ref{C-type-int})  is presented in the monographs by F.D.~Gakhov \cite{Gah} and N.I.~Muskhelishvili \cite{Mus}  under the classical assumptions about the smoothness of the integration curve and the H\"older density of the integral.
In the papers of N.A.~Davydov \cite{Davydov}, V.V.~Salaev \cite{Salaev},
T.S.~Salimov \cite{Salimov}, E.M.~Dyn'kin \cite{Dynkin}, O.F.~Gerus \cite{G77,G96-2-4}, the theory of the Cauchy-type integral and the Cauchy singular integral is developed on an arbitrary rectifiable Jordan curve in extended classes compared to the H\"older class of the integral density,
which are defined, as a rule, in terms of the modulus of continuity of the function $g$.

In the paper of O.F.~Gerus and M.~Shapiro \cite{Ger-Sh-1}, an analog of the Davydov theorem \cite{Davydov} is proved for an appropriate Cauchy-type integral along an arbitrary rectifiable Jordan curve in $\mathbb{R}^2$, which takes values in the algebra of quaternions.
This result is applied in the paper of O.F.~Gerus and M.~Shapiro \cite{Ger-Sh-2} to establish sufficient conditions for the continuous extension to the boundary of a domain of metaharmonic potentials, a partial case of which is the logarithmic double layer potential \eqref{log-pot}.

At the same time, the results mentioned above about the continuous extension of the logarithmic double layer potential
to the boundary of a domain, which are contained in the papers \cite{Plemelj-1911,Radon-46,Kral-1-1964}, are valid for arbitrary continuous functions $g$. It is liked to the fact that the real part of the Schwartz integral, i.e. the Poisson integral, is continuously extended to the boundary of the unit disk for an arbitrary continuous integral density, while the continuous extension of the imaginary part of the Schwartz integral requires additional assumptions about the integral density.

The purpose of this paper is to establish general results about the continuous extension of the real part of the Cauchy-type integral with a real-valued integral density, which are usable for the cases where the classical results of the papers \cite{Plemelj-1911,Radon-46,Kral-1-1964} as well as the corresponding result of the paper \cite{Ger-Sh-2} are not applicable, generally speaking.

\section{Preliminary information}
\label{Preliminary}

In what follows, a closed rectifiable Jordan curve $\gamma$ satisfies the condition (see V.V.~Salaev \cite{Salaev})
\begin{equation}\label{2-4:nerivnist'-teta}
\theta(\varepsilon):=\sup\limits_{\xi\in\gamma} \
\theta_{\xi}(\varepsilon)=O(\varepsilon),\quad
\varepsilon\rightarrow 0\,,
\end{equation}
where $\theta_{\xi}(\varepsilon):={\rm mes}\,\gamma_{\varepsilon}(\xi)$,\,\,
$\gamma_{\varepsilon}(\xi):=\{t\in\gamma :
|t-\xi|\le \varepsilon\}$ and\, ${\rm mes}$\, denotes the linear Lebesgue measure on $\gamma$.
Curves satisfying the condition (\ref{2-4:nerivnist'-teta}) are important in solving various problems
(see, for example, V.V.~Salaev \cite{Salaev}, L.~Ahlfors \cite{Ahlfors2-2-4}, G.~David \cite{David-2-4},
C.~Pommerenke \cite{Pomeren-6-2}, A.~B\"ottcher and Y.I.~Karlovich \cite{Bot-Karl-2-4}).
Such curves are often called {\it regular} (see, for example, \cite{David-2-4}) or {\it Ahlfors-regular}
(see, for example, \cite{Pomeren-6-2}), or {\it Carleson curves} (see, for example, \cite{Bot-Karl-2-4}).

It is well known that the closed  rectifiable Jordan curve $\gamma$ has a tangent at almost all points $t\in\gamma$.
For such a point $t\in\gamma$ we denote by $\vartheta_t$ the angle between the tangent to the curve $\gamma$ at this point and the direction of the real axis. J.~Radon \cite{Radon-46} called $\gamma$ a {\it curve of bounded rotation} if the angle $\vartheta_t$ is a function of bounded variation on $\gamma$.

This implies that for a curve $\gamma$ of bounded rotation, the angle $\vartheta_t$ can have at most a countable set of discontinuity points, and there are one-sided tangents at each point of the curve $\gamma$. Moreover, a curve of bounded rotation can have only a finite set of cusp points and at most a countable set of corner points. At the same time, every curve of bounded rotation satisfies the condition (\ref{2-4:nerivnist'-teta}).

This follows, for example, from the fact that the Cauchy singular integral operator is bounded in Lebesgue spaces on any curve $\gamma$ of bounded rotation (see I.I.~Danilyuk \cite{Dan-m}, I.I.~Daniljuk and V.Ju.~\v{S}elepov \cite{Dan-Shelep}, \`E. G.~Gordadze \cite{Gordadze}), and a necessary condition for this is the condition (\ref{2-4:nerivnist'-teta}) on the curve (see V.A.~Paatashvili and G.A.~Khuskivadze \cite{Paat-Khuskivadze}).

A curve $\gamma$ is called a {\it Lyapunov curve} if the angle $\vartheta_t$ satisfies the H\"older condition:
\[ |\vartheta_{t_1}-\vartheta_{t_2}|\le c\,|t_1-t_2|^{\alpha} \qquad \forall\,t_1,t_2\in \gamma, \]
where $\alpha\in(0,1]$ and the constant\, $c$\, does not depend on $t_1$ and $t_2$.
It is clear that the Lyapunov curve is a smooth curve and also satisfies
the condition (\ref{2-4:nerivnist'-teta}). There are Lyapunov curves that are not the Radon curves of bounded rotation
(see, for example, I.I.~Danilyuk \cite[p.~26]{Dan-m}).

J.~Kr\'al \cite{Kral-1-1964} proved that the logarithmic double layer potential \eqref{log-pot} is extended continuously from the domains
$D^{+}$ and $D^{-}$ to the boundary $\gamma$ for all continuous functions $g : \gamma\rightarrow\mathbb{R}$
if and only if the curve $\gamma$ satisfies the condition
 \begin{equation}\label{um-Kral}
\sup_{\xi\in\gamma}\,\int\limits_0^{2\pi} \mu_{\gamma}(\xi,\phi)\,d\phi<\infty,
\end{equation}
where $\mu_{\gamma}(\xi,\phi)$  is the number of intersection points of the curve $\gamma$ with the ray $\{z=\xi+re^{i\phi} : r>0\}$.
It will be shown below that each curve $\gamma$, which satisfies the condition \eqref{um-Kral},
also satisfies the condition \eqref{2-4:nerivnist'-teta}, i.e, it is an Ahlfors-regular curve.

Note that not every smooth curve satisfies the condition \eqref{um-Kral}.
In particular, an example of such a curve will be given below.

For a function $f \colon E\rightarrow \mathbb{C}$ continuous on a set $E\subset \mathbb{C}$, we shall use its modulus of continuity
\[\omega_{E}(f, \varepsilon):=\sup\limits_{t_1, t_2 \in E\,:\, |t_1-t_2|\le \varepsilon}
\left|f(t_1)-f(t_2)\right|. \]

Quite often, the conditions for a given domain are formulated in terms of a mapping of the unit disk onto this domain (see, for example,
I.I.~Priwalow \cite[\S\S~12--17 of Chapter III]{Priv-6-1}).

Consider a conformal mapping $\sigma : U\rightarrow D^+$ of the unit disk $U$ onto the domain $D^+$. It is well known that the mapping $\sigma$
is continuously extended to the boundary $\partial U$ and defines a homeomorphism between the unit circle $\partial U$ and the curve
$\gamma$.

In the paper  \cite{GrPl-Flaut-2021-3-6}, when solving boundary value problems for monogenic hypercomplex functions associated with a biharmonic equation, the following result on the continuous extension of the logarithmic double layer potential (\ref{log-pot}) has actually been established, although it was not formulated as a theorem:

\vskip 1mm

\begin{theorem}\label{theor-0}
{\it Let a conformal mapping $\sigma : U\rightarrow D^+$ of the unit disk $U$ onto the domain $D^+$ have the nonvanishing continuous contour derivative $\sigma\,'$ on the circle $\partial U$, and let its modulus of continuity satisfy the Dini
condition
\begin{equation}\label{um-Dini-conf-vid}
\int\limits_{0}^{1}\frac{\omega_{\partial U}(\sigma\,',\eta)}{\eta}\,d\eta<\infty.
\end{equation}
Then, for each continuous function $g \colon \gamma\rightarrow \mathbb{R}$, the integral \eqref{log-pot} is continuously extended from the domain
$D^+$ to the boundary $\gamma$.}
\end{theorem}

\vskip 1mm

Theorem \ref{theor-0} generalizes the corresponding result of the classical theory of the logarithmic double layer potential on the Lyapunov curves (see, for example, J.~Plemelj \cite{Plemelj-1911}), 
because in the case where $\gamma$ is a Lyapunov curve,
the condition (\ref{um-Dini-conf-vid}) is satisfied owing to the Kellogg theorem (see, for example, G.M.~Goluzin \cite{Gol-6-3}).
The condition (\ref{um-Dini-conf-vid}) is also satisfied in the more general case where the modulus of continuity of
the angle\, $\vartheta_t$\, satisfies the condition
\begin{equation}\label{um-Dini-ln-dotychna}
\int\limits_{0}^{1}\frac{\omega_{\gamma}(\vartheta_t,\eta)}{\eta} \
\ln\frac{2}{\eta}\,d\eta<\infty\,.
 \end{equation}
It follows from the estimate of the modulus of continuity of the function $\sigma\,'$
presented in Theorem 2 in the paper of J.L.~Heronimus \cite{Grnm-6-3}\, (see also S.E.~Warschawski \cite{War2-6-3}).

If the modulus of continuity of the function\, $g : \gamma\rightarrow\mathbb{R}$\, satisfies the Dini
condition
 \begin{equation}\label{2-4:Dini}
\int\limits_0^1\frac{\omega_{\gamma}(g,\eta)}{\eta}\,d\eta<\infty,
\end{equation}
then the reduced singular Cauchy integral
\begin{equation}\label{zved-siC}
\int\limits_{\gamma}\frac{g(t)-g(\xi)}{t-\xi}\,dt:=\lim\limits_{\delta\to 0^+}\,\int\limits_{\gamma\setminus\gamma_{\delta}(\xi)}
 \frac{g(t)-g(\xi)}{t-\xi}\,dt \qquad \forall\,\xi\in\gamma
 \end{equation}
exists (see O.F.~Gerus \cite{G77}, and also V.V.~Salaev \cite{Salaev}, where the Dini
condition of the form (\ref{2-4:Dini}) is given in terms of the regularized modulus of continuity using
the Stechkin construction).

From this and from the result of N.A.~Davydov \cite{Davydov}, it follows that the Cauchy-type integral (\ref{C-type-int}) has
the limiting values $\widetilde{g}\,^{\pm}(\xi)$\, at every point\, $\xi\in\gamma$\, from the domains\, $D^{\pm}$, which are expressed by the Sokhotski--Plemelj formulas:
\begin{equation}\label{2-4:Dav-Sokh+}
 \widetilde{g}\,^{+}(\xi)=g(\xi)+\frac{1}{2\pi i}\,\int\limits_{\gamma}\frac{g(t)-g(\xi)}{t-\xi}\,dt\,,
\end{equation}
\begin{equation}\label{2-4:Dav-Sokh-}
 \widetilde{g}\,^{-}(\xi)=\frac{1}{2\pi i}\,\int\limits_{\gamma}\frac{g(t)-g(\xi)}{t-\xi}\,dt\,.
\end{equation}

Let $d:=\max\limits_{t_1,t_2\in\gamma} |t_1-t_2|$ be the diameter of the curve $\gamma$.

To single out the real part of the integral (\ref{zved-siC}), we define a branch $\arg(z-\xi)$ continuous on $\gamma\setminus\{\xi\}$
of the multivalued function ${\rm Arg}\,(z-\xi)$ in the following way.
For each positive $\delta<d/2$, we select that connected component $\gamma_{\xi,\delta}$ of the set $\gamma_{\delta}(\xi)$ which contains the point $\xi$, and we take such a point $\xi_1\in\gamma_{\xi,\delta}$ at which there is a tangent to $\gamma$ and which does not precede the point $\xi$ under the given orientation of the curve $\gamma$. It is obvious that in the case where there is a tangent to $\gamma$ at the point $\xi$, we can set $\xi_1=\xi$.
Let us cut the complex plane along the curve $\Gamma_{\xi,\delta}:=\gamma[\xi,\xi_1]\cup\Gamma[\xi_1,\infty]$, where $\gamma[\xi,\xi_1]$ is the arc of $\gamma$ with the initial point $\xi$ and the end point $\xi_1$, and $\Gamma[\xi_1,\infty]$ is a smooth curve that connects the points $\xi_1$ and $\infty$ and lies completely (except for its ends $\xi_1$ and $\infty$) in the domain $D^-$.
Now, let us single out a branch $\arg_{\delta}(z-\xi)$ of the multivalued function ${\rm Arg}\,(z-\xi)$, which is continuous outside the cut
$\Gamma_{\xi,\delta}$ with the normalization condition $\arg_{\delta}(z_0-\xi)=\phi_0$, where $z_0\in D^+$ and $\phi_0$ is one of the values of the function ${\rm Arg}\,(z-\xi)$ at $z=z_0$. We shall use the fixed values $z_0$ and $\phi_0$ for all positive $\delta<d/2$.
As a result, we have the obvious equality
\[\arg_{\delta_1}(t-\xi)=\arg_{\delta}(t-\xi) \qquad \forall\,\delta_1,\delta : 0<\delta_1<\delta<d/2 \quad \forall\,t\in\gamma\setminus\gamma_{\delta}(\xi)\]
that implies the existence of the following limit:
\[\arg(t-\xi):=\lim\limits_{\delta\to 0^+}\,\arg_{\delta}(t-\xi) \qquad \forall\,t\in\gamma\setminus\{\xi\}\,.\]

Thus, under the assumption that the function\, $g : \gamma\rightarrow\mathbb{R}$\, satisfies the condition (\ref{2-4:Dini}), from the equality (\ref{zved-siC}) we get the equality
\[{\rm Re}\,\bigg(\frac{1}{2\pi i}\,\int\limits_{\gamma}\frac{g(t)-g(\xi)}{t-\xi}\,dt\bigg)=\frac{1}{2\pi}\,\lim\limits_{\delta\to 0^+}\,\int\limits_{\gamma\setminus\gamma_{\delta}(\xi)}
 \big(g(t)-g(\xi)\big)\,d\arg(t-\xi) \qquad \forall\,\xi\in\gamma\,,\]
where on the right-hand side of the equality the integral is the Stieltjes integral.

Let us accept by definition
\begin{equation}\label{Re-siC}
\int\limits_{\gamma}
 \big(g(t)-g(\xi)\big)\,d\arg(t-\xi):=\lim\limits_{\delta\to 0^+}\,\int\limits_{\gamma\setminus\gamma_{\delta}(\xi)}
 \big(g(t)-g(\xi)\big)\,d\arg(t-\xi) \qquad \forall\,\xi\in\gamma\,.
  \end{equation}

Finally, denoting
\[\big({\rm Re}\,\widetilde{g}\big)^{\pm}(\xi):=\lim\limits_{z\to\xi,\,z\in D^{\pm}}\,{\rm Re}\,\widetilde{g}(z)\qquad \forall\,\xi\in\gamma,\]
as a consequence of the formulas (\ref{2-4:Dav-Sokh+}) and (\ref{2-4:Dav-Sokh-}), for all $\xi\in\gamma$, we obtain the equalities
\begin{equation}\label{double-pot+}
 \big({\rm Re}\,\widetilde{g}\big)^{+}(\xi)=g(\xi)+\frac{1}{2\pi}\,\int\limits_{\gamma}
 \big(g(t)-g(\xi)\big)\,d\arg(t-\xi)\,,
\end{equation}
\begin{equation}\label{double-pot-}
 \big({\rm Re}\,\widetilde{g}\big)^{-}(\xi)=\frac{1}{2\pi}\,\int\limits_{\gamma}
 \big(g(t)-g(\xi)\big)\,d\arg(t-\xi)\,.
\end{equation}

Below, we investigate the fulfillment of the equalities (\ref{double-pot+}) and (\ref{double-pot-}), not assuming, generally speaking, neither the fulfillment of the condition (\ref{2-4:Dini}) for the function\, $g : \gamma\rightarrow\mathbb{R}$\, nor the fulfillment of the condition \eqref{um-Kral} for the curve $\gamma$.

\section{A necessary and sufficient condition for the continuous extension of the real part of the Cauchy-type integral to the boundary of the domain bounded by an Ahlfors-regular curve }
\label{Criterion}

The following statement is true:

\vskip 1mm

\begin{theorem}\label{theor-crit}
{\it Let a closed Jordan curve $\gamma$ be Ahlfors-regular
and let a function
$g \colon \gamma\rightarrow \mathbb{R}$ be continuous on $\gamma$.
The function ${\rm Re}\,\widetilde{g}(z)$ is continuously extended to the boundary $\gamma$ from the domain $D^+$ or $D^-$ if and only if the following condition is satisfied:
\begin{equation}\label{nidu-nepr}
\sup\limits_{\xi\in\gamma}\,\sup\limits_{\delta\in(0,\varepsilon)}\,\bigg|\,\int\limits_{\gamma_{\varepsilon}(\xi)\setminus\gamma_{\delta}(\xi)}
 \big(g(t)-g(\xi)\big)\,d\arg(t-\xi)\,\bigg|\to 0\,, \qquad \varepsilon\to 0\,.
\end{equation}
In the case where the condition \eqref{nidu-nepr} is satisfied, the limiting values $\big({\rm Re}\,\widetilde{g}\big)^{\pm}(\xi)$ are represented by the formulas \eqref{double-pot+} and \eqref{double-pot-} for all $\xi\in\gamma$.}
\end{theorem}

\vskip 1mm

\textit{\textbf{Proof.}} {\it Sufficiency.} Obviously, when the condition (\ref{nidu-nepr}) is satisfied, the limit exists in the equality (\ref{Re-siC}).

Let us prove the equality (\ref{double-pot+}). Let $\xi\in\gamma$, $z\in D^+$ and $\varepsilon:=|z-\xi|<d/8$.
Denote $\varepsilon_1:=\min\limits_{t\in\gamma}|t-z|$.
Let us choose a point $\xi_z\in\gamma$ closest to the point $z$.

We use the representation of the difference
  \begin{multline*}
{\rm Re}\,\widetilde{g}(z)-g(\xi)-\frac{1}{2\pi}\,\int\limits_{\gamma}
 \big(g(t)-g(\xi)\big)\,d\arg(t-\xi)\\[2mm]
  ={\rm Re}\,\bigg(\frac{1}{2\pi i}\,\int\limits_{\gamma}\frac{g(t)-g(\xi_z)}{t-z}\,dt\bigg)-\frac{1}{2\pi}\,\int\limits_{\gamma}
 \big(g(t)-g(\xi_z)\big)\,d\arg(t-\xi_z)+g(\xi_z)-g(\xi)\\[2mm]
 +\frac{1}{2\pi}\,\int\limits_{\gamma}
 \big(g(t)-g(\xi_z)\big)\,d\arg(t-\xi_z)-\frac{1}{2\pi}\,\int\limits_{\gamma}
 \big(g(t)-g(\xi)\big)\,d\arg(t-\xi)\,.
\end{multline*}

Consider the difference
 \begin{multline*}
 {\rm Re}\,\bigg(\frac{1}{2\pi i}\,\int\limits_{\gamma}\frac{g(t)-g(\xi_z)}{t-z}\,dt\bigg)-\frac{1}{2\pi}\,\int\limits_{\gamma}
 \big(g(t)-g(\xi_z)\big)\,d\arg(t-\xi_z)
 ={\rm Re}\,\bigg(\frac{1}{2\pi i}\,\int\limits_{\gamma_{2\varepsilon_1}(\xi_z)}\frac{g(t)-g(\xi_z)}{t-z}\,dt\bigg)\\
 - \frac{1}{2\pi}\,\int\limits_{_{2\varepsilon_1}(\xi_z)}\big(g(t)-g(\xi_z)\big)\,d\arg(t-\xi_z)
 +{\rm Re}\,\bigg(\frac{z-\xi_z}{2\pi i}\,\int\limits_{\gamma\setminus\gamma_{2\varepsilon_1}(\xi_z)}\frac{g(t)-g(\xi_z)}{(t-z)(t-\xi_z)}\,dt\bigg)
=:I_1-I_2+I_3\,.
\end{multline*}

Taking into account the condition (\ref{2-4:nerivnist'-teta}), we obtain the relation
\[ |I_1|\le \frac{1}{2\pi}\,\int\limits_{\gamma_{2\varepsilon_1}(\xi_z)}\frac{|g(t)-g(\xi_z)|}{|t-z|}\,|dt|\le
 \frac{\omega_{\gamma}(g,2\varepsilon_1)}{2\pi\varepsilon_1}\,\theta_{\xi_z}(2\varepsilon_1)\le c\,\omega_{\gamma}(g,2\varepsilon_1)
 \rightarrow 0, \qquad \varepsilon_1\rightarrow 0, \]
where the constant\, $c$\, depends only on the curve $\gamma$.

The condition \eqref{nidu-nepr} implies the relation
\[ |I_2|\rightarrow 0, \qquad \varepsilon_1\rightarrow 0. \]

To estimate the integral $I_3$, we use Proposition 7.2 in \cite{Pla-Shpak-mono}
(see also the proof of Theorem 1 in the paper of O.F.~Gerus \cite{G78}) and the condition
(\ref{2-4:nerivnist'-teta}) so that we have
\begin{multline*}
|I_3|\le \frac{|z-\xi_z|}{\pi}\,\int\limits_{\gamma\setminus\gamma_{2\varepsilon_1}(\xi_z)}\frac{|g(t)-g(\xi_z)|}{|t-\xi_z|^2}\,|dt|\le
\frac{\varepsilon_1}{\pi}\int\limits_{[2\varepsilon_1,d]}\frac{\omega_{\gamma}(g,\eta)}{\eta^2}\,d\theta_{\xi_z}(\eta)\\
\leq \frac{2\varepsilon_1}{3\pi}\,\int\limits_{\varepsilon_1}^{d}\frac{\theta_{\xi_z}(2\eta)\omega_{\gamma}(g,2\eta)}{\eta^3}\,d\eta\leq
c\,\varepsilon_1\int\limits_{\varepsilon_1}^{2d}\frac{\omega_{\gamma}(g,\eta)}{\eta^2}\,d\eta \rightarrow 0, \qquad \varepsilon_1\rightarrow 0,
\end{multline*}
where the constant\, $c$\, depends only on the curve $\gamma$.

Now, consider the difference
 \begin{multline*}
  \frac{1}{2\pi}\,\int\limits_{\gamma}
 \big(g(t)-g(\xi_z)\big)\,d\arg(t-\xi_z)-\frac{1}{2\pi}\,\int\limits_{\gamma}
 \big(g(t)-g(\xi)\big)\,d\arg(t-\xi)\\
 =\frac{1}{2\pi}\,\int\limits_{\gamma_{\varepsilon}(\xi_z)}
 \big(g(t)-g(\xi_z)\big)\,d\arg(t-\xi_z)+{\rm Re}\,\bigg(\frac{1}{2\pi i}\,\int\limits_{\gamma_{4\varepsilon}(\xi)\setminus\gamma_{\varepsilon}(\xi_z)}\frac{g(t)-g(\xi_z)}{t-\xi_z}\,dt\bigg)
  -\frac{1}{2\pi}\,\int\limits_{\gamma_{4\varepsilon}(\xi)}\big(g(t)-g(\xi)\big)\,d\arg(t-\xi)\\
 +{\rm Re}\,\bigg(\frac{\xi_z-\xi}{2\pi i}\, \int\limits_{\gamma\setminus\gamma_{4\varepsilon}(\xi)}\frac{g(t)-g(\xi)}{(t-\xi)(t-\xi_z)}\,dt\bigg)
 +{\rm Re}\,\bigg(\frac{g(\xi)-g(\xi_z)}{2\pi i}\, \int\limits_{\gamma\setminus\gamma_{4\varepsilon}(\xi)}\frac{dt}{t-\xi_z}\bigg)
 =:J_1+J_2-J_3+J_4+J_5\,.
\end{multline*}

The condition  \eqref{nidu-nepr} implies the relations
\[ |J_1|\rightarrow 0 \quad \mbox{and} \quad |J_3|\rightarrow 0, \qquad \varepsilon\rightarrow 0. \]

The integrals $J_2$ and $J_4$ are estimated similarly to the integrals $I_1$ and $I_3$, respectively. As a result, we have the relations
\[ |J_2|\rightarrow 0 \quad \mbox{and} \quad |J_4|\rightarrow 0, \qquad \varepsilon\rightarrow 0. \]

In addition, the following relations are fulfilled (see the proof of Theorem 1 in the paper of V.V.~Salaev \cite{Salaev}):
\[|J_5|\le \frac{|g(\xi)-g(\xi_z)|}{2\pi}\,\bigg|\, \int\limits_{\gamma\setminus\gamma_{4\varepsilon}(\xi)}\frac{dt}{t-\xi_z}\,\bigg|\le
2\,\omega_{\gamma}(g,2\varepsilon)\rightarrow 0, \qquad \varepsilon\rightarrow 0.\]

An obvious consequence of the given relations is the equality \eqref{double-pot+}. The equality \eqref{double-pot-} is similarly established.

{\it Necessity.} Since the curve $\gamma$ satisfies the condition \eqref{2-4:nerivnist'-teta}, the singular Cauchy integral operator is bounded in the Lebesgue spaces $L_p$ for $p>1$ on $\gamma$ (see Theorem 1 in the paper of G.~David \cite{David-2-4}). At the same time, the Cauchy-type integral \eqref{C-type-int} belongs to the Smirnov classes $E_p$ (see, for example,
I.I.~Priwalow \cite{Priv-6-1}) for $p>1$ in the domains $D^+$ and $D^-$. In addition, its angular boundary values
$\widetilde{g}\,^{\pm}_{\rm ang}(\xi)$\, from $D^{\pm}$ exist for almost all points $\xi\in\gamma$, and the following equality holds almost everywhere on $\gamma$:
\[g(\xi)=\widetilde{g}\,^{+}_{\rm ang}(\xi)-\widetilde{g}\,^{-}_{\rm ang}(\xi). \]

We denote by $\widetilde{g}\,^{\pm}$ the function that is defined by the equality \eqref{C-type-int} in the domain $D^{\pm}$ and is extended almost everywhere on $\gamma$ by means of the values $\widetilde{g}\,^{\pm}_{\rm ang}$.
We denote also the real part of this function by ${\rm Re}\,\widetilde{g}\,^{\pm}$.

Note that the values of the functions ${\rm Re}\,\widetilde{g}\,^{+}$ and ${\rm Re}\,\widetilde{g}\,^{-}$  are expressed by equalities of the form \eqref{double-pot+} and \eqref{double-pot-} for almost all points $\xi\in\gamma$.
It obviously follows that in the case where the function ${\rm Re}\,\widetilde{g}(z)$ is continuously extended to the boundary $\gamma$ from one of the domains $D^+$ or $D^-$, this function is also continuously extended to $\gamma$ from the other domain.

For $\xi\in\gamma$ and $0<\delta<\varepsilon<d$, consider the open sets
$D^{\pm}_{\delta,\varepsilon}(\xi):=\{z\in D^{\pm} : \delta<|z-\xi|<\varepsilon\}$ and their boundaries $\partial D^{\pm}_{\delta,\varepsilon}(\xi)$, the orientation of which is induced by the orientation of 
$\gamma$. Denote $\Gamma^{\pm}_{\delta}:=\{t\in\partial D^{\pm}_{\delta,\varepsilon}(\xi)\setminus\gamma : |z-\xi|=\delta\}$, $\Gamma^{\pm}_{\varepsilon}:=\{t\in\partial D^{\pm}_{\delta,\varepsilon}(\xi)\setminus\gamma : |z-\xi|=\varepsilon\}$.

We have the equalities
\begin{multline*}
\int\limits_{\gamma_{\varepsilon}(\xi)\setminus\gamma_{\delta}(\xi)}
 \big(g(t)-g(\xi)\big)\,d\arg(t-\xi)={\rm Im}\,\bigg(\,\int\limits_{\gamma_{\varepsilon}(\xi)\setminus\gamma_{\delta}(\xi)}
 \frac{g(t)-g(\xi)}{t-\xi}\,dt\bigg)\\[2mm]
 ={\rm Im}\,\bigg(\,\int\limits_{\gamma_{\varepsilon}(\xi)\setminus\gamma_{\delta}(\xi)}
 \frac{\widetilde{g}\,^{+}(t)-\widetilde{g}\,^{-}(t)-{\rm Re}\,\widetilde{g}\,^{+}(\xi)+{\rm Re}\,\widetilde{g}\,^{-}(\xi)}{t-\xi}\,dt\bigg)\\[2mm]
={\rm Im}\,\bigg(\,\int\limits_{\partial D^{+}_{\delta,\varepsilon}(\xi)}
 \frac{\widetilde{g}\,^{+}(t)-{\rm Re}\,\widetilde{g}\,^{+}(\xi)}{t-\xi}\,dt
 -\int\limits_{\Gamma^{+}_{\delta}} \frac{\widetilde{g}\,^{+}(t)- {\rm Re}\,\widetilde{g}\,^{+}(\xi)}{t-\xi}\,dt-\int\limits_{\Gamma^{+}_{\varepsilon}}
 \frac{\widetilde{g}\,^{+}(t)-{\rm Re}\,\widetilde{g}\,^{+}(\xi)}{t-\xi}\,dt\bigg) 
\end{multline*}
\begin{equation*}
 -{\rm Im}\,\bigg(\,\int\limits_{\partial D^{-}_{\delta,\varepsilon}(\xi)}
 \frac{\widetilde{g}\,^{-}(t)-{\rm Re}\,\widetilde{g}\,^{-}(\xi)}{t-\xi}\,dt
 -\int\limits_{\Gamma^{-}_{\delta}}
 \frac{\widetilde{g}\,^{-}(t)-{\rm Re}\,\widetilde{g}\,^{-}(\xi)}{t-\xi}\,dt-\int\limits_{\Gamma^{-}_{\varepsilon}}
 \frac{\widetilde{g}\,^{-}(t)-{\rm Re}\,\widetilde{g}\,^{-}(\xi)}{t-\xi}\,dt\bigg).
\end{equation*}

Further, taking into account that the integrals of functions from the Smirnov classes along the closed curves $\partial D^{\pm}_{\delta,\varepsilon}(\xi)$ are equal to zero, we have
 \begin{multline*}
\int\limits_{\gamma_{\varepsilon}(\xi)\setminus\gamma_{\delta}(\xi)}
 \big(g(t)-g(\xi)\big)\,d\arg(t-\xi)
 =-{\rm Im}\,\bigg(\,\int\limits_{\Gamma^{+}_{\delta}} \frac{\widetilde{g}\,^{+}(t)- {\rm Re}\,\widetilde{g}\,^{+}(\xi)}{t-\xi}\,dt+\int\limits_{\Gamma^{+}_{\varepsilon}}
 \frac{\widetilde{g}\,^{+}(t)-{\rm Re}\,\widetilde{g}\,^{+}(\xi)}{t-\xi}\,dt\bigg)\\
 +{\rm Im}\,\bigg(\,\int\limits_{\Gamma^{-}_{\delta}}\frac{\widetilde{g}\,^{-}(t)-{\rm Re}\,\widetilde{g}\,^{-}(\xi)}{t-\xi}\,dt+ \int\limits_{\Gamma^{-}_{\varepsilon}}\frac{\widetilde{g}\,^{-}(t)-{\rm Re}\,\widetilde{g}\,^{-}(\xi)}{t-\xi}\,dt\bigg)\\
=-{\rm Im}\,\bigg(\,\int\limits_{\Gamma^{+}_{\delta}} \frac{{\rm Re}\,\widetilde{g}\,^{+}(t)-{\rm Re}\,\widetilde{g}\,^{+}(\xi)}{t-\xi}\,dt+ \int\limits_{\Gamma^{+}_{\varepsilon}}
 \frac{{\rm Re}\,\widetilde{g}\,^{+}(t)-{\rm Re}\,\widetilde{g}\,^{+}(\xi)}{t-\xi}\,dt\bigg)\\
 +{\rm Im}\,\bigg(\,\int\limits_{\Gamma^{-}_{\delta}}\frac{{\rm Re}\,\widetilde{g}\,^{-}(t)-{\rm Re}\,\widetilde{g}\,^{-}(\xi)}{t-\xi}\,dt+ \int\limits_{\Gamma^{-}_{\varepsilon}}\frac{{\rm Re}\,\widetilde{g}\,^{-}(t)-{\rm Re}\,\widetilde{g}\,^{-}(\xi)}{t-\xi}\,dt\bigg).
   \end{multline*}

Since the function ${\rm Re}\,\widetilde{g}(z)$ is continuously extended to the boundary $\gamma$ from $D^{\pm}$ and
vanishes at infinity, the function ${\rm Re}\,\widetilde{g}\,^{\pm}$ is uniformly continuous in the closure
$\overline{D^{\pm}}$ of the domain $D^{\pm}$. Therefore, we obtain the estimations
\begin{multline*}
\bigg|\,\int\limits_{\gamma_{\varepsilon}(\xi)\setminus\gamma_{\delta}(\xi)}
 \big(g(t)-g(\xi)\big)\,d\arg(t-\xi)\bigg|
\le\int\limits_{\Gamma^{+}_{\delta}} \frac{\big|{\rm Re}\,\widetilde{g}\,^{+}(t)-{\rm Re}\,\widetilde{g}\,^{+}(\xi)\big|}{|t-\xi|}\,|dt|+ \int\limits_{\Gamma^{+}_{\varepsilon}}
 \frac{\big|{\rm Re}\,\widetilde{g}\,^{+}(t)-{\rm Re}\,\widetilde{g}\,^{+}(\xi)\big|}{|t-\xi|}\,|dt|\\
 +\int\limits_{\Gamma^{-}_{\delta}}\frac{\big|{\rm Re}\,\widetilde{g}\,^{-}(t)-{\rm Re}\,\widetilde{g}\,^{-}(\xi)\big|}{|t-\xi|}\,|dt|+ \int\limits_{\Gamma^{-}_{\varepsilon}}\frac{\big|{\rm Re}\,\widetilde{g}\,^{-}(t)-{\rm Re}\,\widetilde{g}\,^{-}(\xi)\big|}{|t-\xi|}\,|dt|
 \le 4\pi\,\omega\,_{\overline{D^+}}\,({\rm Re}\,\widetilde{g}\,^+,\varepsilon)+4\pi\,\omega\,_{\overline{D^-}}\,({\rm Re}\,\widetilde{g}\,^-,\varepsilon),
    \end{multline*}
which imply the condition \eqref{nidu-nepr}.
\hfill $\Box$

Theorem \ref{theor-crit} is similar in a certain sense to the corresponding theorem for the Cauchy-type integral, which is proved by A.O.~Tokov
\cite{Tokov}.

\section{Some properties of Ahlfors-regular curves}
\label{Properties of Ahlfors-regular}

Note that for each $\xi\in\gamma$ and each $\delta>0$, the function $\arg(t-\xi)$ has a bounded variation on the set $\gamma\setminus\gamma_{\delta}(\xi)$. However, in general, the function $\arg(t-\xi)$ can be a function of unbounded variation on $\gamma$, because, in particular, it can be unbounded in a neighborhood of the point $\xi$.

Consider the class of curves  $\gamma$, for which the function $\arg(t-\xi)$ has a bounded variation  $V_{\gamma}[\arg(t-\xi)]$ on $\gamma\setminus\{\xi\}$ for all  $\xi\in\gamma$ and, moreover, satisfies the condition
\begin{equation}\label{obm-var-arg}
\sup\limits_{\xi\in\gamma} V_{\gamma}[\arg(t-\xi)]<\infty\,.
\end{equation}
It is obvious that a curve satisfying the condition \eqref{obm-var-arg} has one-sided tangents at each point $\xi\in\gamma$.

Note that the condition \eqref{obm-var-arg} is equivalent to the condition \eqref{um-Kral}, which follows from the Banach indicatrix theorem (see J.~Kr\'al \cite[Lemma 1.2]{Kral-1-1964}). Thus, the class of curves satisfying the condition \eqref{obm-var-arg} includes curves from the corresponding classical results of J.~Plemelj \cite{Plemelj-1911} and J.~Radon \cite{Radon-46} and from Theorem \ref{theor-0}.

Curves satisfying the condition \eqref{obm-var-arg} will be called the {\it Kr\'al curves}.

\vskip 1mm

\begin{proposition}\label{prop-0}
{\it Every Kr\'al curve is an Ahlfors-regular curve.}
\end{proposition}

\vskip 1mm

\textit{\textbf{Proof.}} Let $\gamma$ be a Kr\'al curve and $\xi\in\gamma$. Let us first show that for an arbitrary $\varepsilon>0$,
the variation of the function $|t-\xi|$ on the set $\gamma_{\varepsilon}(\xi)$ satisfies the inequality
\begin{equation}\label{reg-mod}
V_{\gamma_{\varepsilon}(\xi)}[|t-\xi|] 
\le c\,\varepsilon,
\end{equation}
where the constant\, $c$\, does not depend on $\xi$ and $\varepsilon$.

In the case where we consider a fixed point $z\in\mathbb{C}\setminus\gamma$ and a variable $t\in\gamma$, we understand $\arg(t-z)$ as an arbitrary branch of the multivalued function  ${\rm Arg}\,(t-z)$.

We use the representation $\gamma_{\varepsilon}(\xi)=\gamma_1\cup\gamma_2\cup\gamma_3$, where
\[\gamma_1:=\gamma_{\varepsilon}(\xi)\cap\big\{t=\xi+r\,e^{i\phi} : r>0,\, \phi\in(-\pi/4,\pi/4)\cup(3\pi/4,5\pi/4)\big\},  \]
\[\gamma_2:=\gamma_{\varepsilon}(\xi)\cap\big\{t=\xi+r\,e^{i\phi} : r>0,\, \phi\in(-3\pi/4,-\pi/4)\cup(\pi/4,3\pi/4)\big\},  \]
\[\gamma_3:=\gamma_{\varepsilon}(\xi)\cap\big\{t=\xi+r\,e^{i\phi} : r\ge 0,\, \phi\in\{-\pi/4,\pi/4,-3\pi/4,3\pi/4\}\big\}.  \]

For the variation of the function $|t-\xi|$ on the set $\gamma_1$, the following inequality holds
(see J.~Kr\'al \cite[Theorem 2.10]{Kral-0-1964}):
\[\frac{V_{\gamma_1}[|t-\xi|]}{\varepsilon}\le c_0\Big(V_{\gamma_1}[\arg(t-\xi)]+V_{\gamma_1}[\arg(t-\xi-\varepsilon)] \Big), \]
where  $c_0=6/\sin^2 (\pi/4)=12$.
Moreover, since the curve $\gamma$ satisfies the condition  \eqref{obm-var-arg}, by virtue of Theorem 1.11 from the paper
J.~Kr\'al \cite{Kral-1-1964}, the following condition is also satisfied:
\[\sup\limits_{z\in\mathbb{C}} V_{\gamma}[\arg(t-z)]<\infty\,.\]
As a result, under the condition \eqref{obm-var-arg} for the curve $\gamma$, we obtain the inequality
\begin{equation}\label{reg-mod-1}
V_{\gamma_1}[|t-\xi|] \le c_1\,\varepsilon,
\end{equation}
where the constant\, $c_1$\, does not depend on $\xi$ and $\varepsilon$.

In a similar way, for the variation of the function $|t-\xi|$ on the set $\gamma_2$,  we obtain the inequalities
\begin{equation}\label{reg-mod-2}
V_{\gamma_2}[|t-\xi|]\le c_0\Big(V_{\gamma_2}[\arg(t-\xi)]+V_{\gamma_2}[\arg(t-\xi-i\varepsilon)] \Big)\,\varepsilon \le c_1\,\varepsilon.
\end{equation}

In addition, it is obvious that
\begin{equation}\label{reg-mod-3}
V_{\gamma_3}[|t-\xi|] \le 4\,\varepsilon.
\end{equation}

The inequalities \eqref{reg-mod-1} -- \eqref{reg-mod-3} imply the inequality \eqref{reg-mod}, where $c=2c_1 +4$.

Now, taking into account the condition \eqref{obm-var-arg} and the inequality \eqref{reg-mod}, for arbitrary  $\xi\in\gamma$ and $\varepsilon>0$, we obtain the relations
\begin{equation*}
\theta_{\xi}(\varepsilon)=\,\int\limits_{\gamma_{\varepsilon}(\xi)}|dt|\le \int\limits_{\gamma_{\varepsilon}(\xi)}\big|d|t-\xi|\big|+
\int\limits_{\gamma_{\varepsilon}(\xi)}|t-\xi||d\arg(t-\xi)|
\le V_{\gamma_{\varepsilon}(\xi)}[|t-\xi|]+\varepsilon\,V_{\gamma}[\arg(t-\xi)]\le c\,\varepsilon,
\end{equation*}
where the constant\, $c$\, does not depend on $\xi$ and $\varepsilon$. Thus, the curve  $\gamma$ satisfies the condition \eqref{2-4:nerivnist'-teta}, i.e., it is an Ahlfors-regular curve.
\hfill $\Box$

Among the Kr\'al curves there are curves that are not the Radon curves of bounded rotation, as the following example shows:

\vskip 1mm

\begin{example}\label{ex-1}
Consider the curve
\begin{multline*}
\gamma=\Big\{z=e^{i\phi} : \phi\in[0,\pi]\Big\}\cup[-1,0]\cup\bigcup_{n=1}^{\infty}[2^{-2n+1},2^{-2n+2}]\\
\cup\bigcup_{n=1}^{\infty} \Big\{z=2^{-n}e^{i\phi} : \phi\in[0,2^{-n}]\Big\}\cup \bigcup_{n=1}^{\infty} \Big\{z=re^{ir} : r\in[2^{-2n},2^{-2n+1}]\Big\}.
\end{multline*}

It is clear that $V_{\gamma}[\vartheta_t]=\infty$, but at the same time, the following relations are fulfilled:
\[ V_{\gamma}[\arg t]\le \pi+\sum\limits_{n=1}^{\infty}2^{-n}+\frac{1}{2}=\pi+\frac{3}{2}\,, \]
\begin{multline*}
 V_{\gamma}[\arg(t-\xi)]=V_{\gamma_{|\xi|/2}(0)}[\arg(t-\xi)]+V_{\gamma_{|\xi|}(\xi)\setminus\gamma_{|\xi|/2}(0)}[\arg(t-\xi)]+
V_{\gamma\setminus\gamma_{|\xi|}(\xi)\setminus\gamma_{|\xi|/2}(0)}[\arg(t-\xi)]\\[3mm]
\le V_{\gamma_{|\xi|/2}(0)}[\arg t] +2\pi+2V_{\gamma\setminus\gamma_{|\xi|}(\xi)\setminus\gamma_{|\xi|/2}(0)}[\arg t]\le 2V_{\gamma}[\arg t]+2\pi \qquad \forall\,\xi\in\gamma : 0<|\xi|<1,\\
\end{multline*}
\[V_{\gamma}[\arg(t-\xi)]\le V_{\gamma_{1/2}(\xi)}[\arg(t-\xi)]+\int\limits_{\gamma\setminus\gamma_{1/2}(\xi)}\frac{|dt|}{|t-\xi|}\le
\pi+ 2\,{\rm mes}\,\gamma \qquad \forall\,\xi\in\gamma : |\xi|=1. \]

Thus, $\gamma$ is a Kr\'al curve that is not the Radon curve of bounded rotation.
\end{example}

\vskip 2mm

For points $\xi_1,\xi_2\in\gamma$, we denote by $\gamma[\xi_1,\xi_2]$ the arc of the curve $\gamma$ with the initial point $\xi_1$ and the end point $\xi_2$ at the orientation of this arc, which is induced by the orientation of the curve $\gamma$.

The following statement defines a class of smooth curves satisfying the condition \eqref{obm-var-arg}:

\vskip 2mm

\begin{proposition}\label{prop-3}
{\it If for a closed smooth Jordan curve  $\gamma$ the angle $\vartheta_t$ satisfies the condition
\begin{equation}\label{um-Dini-dotychna}
 \int\limits_{0}^{1}\frac{\omega_{\gamma}(\vartheta_t,\eta)}{\eta}\,d\eta<\infty,
\end{equation}
then $\gamma$ is a Kr\'al curve.}
\end{proposition}

\vskip 1mm

\textit{\textbf{Proof.}} Let $\xi\in\gamma$. It is known (see, for example, N.I.~Muskhelishvili \cite{Mus}) that there exists $r_0>0$, which does not depend on $\xi$,
such that each circle of radius $r\le r_0$ centered at the point $\xi$ intersects $\gamma$ in only two points.

We denote by $t_-$ and $t_+$ the points of intersection of the circle $\{z\in\mathbb{C} : |z-\xi|=r_0\}$ and the curve $\gamma$, and with the given orientation $\gamma$ the point $t_-$ precedes the point $\xi$ and the point $t_+$ follows it. Considering one of the arcs either  $\gamma[t_-,\xi]$ or $\gamma[\xi,t_+]$, we will denote it $\widetilde{\gamma}$.

The arc $\widetilde{\gamma}$ allows the parameterization \,\,
$t=\xi+r\,e^{i(\phi(r)+\phi_0)}$, \,\, $r\in[0,r_0]$, \,\,
where $\phi_0$ is a real constant and $\phi(r)\to 0$ as $r\to 0$. Denote \,\, $\widetilde{x}(r):=r\cos \phi(r)$, \,\, $\widetilde{y}(r):=r\sin \phi(r)$.

For all $t\in\widetilde{\gamma}\setminus\{\xi\}$, the following equalities hold:
\begin{multline*}
d\arg(t-\xi)=d\phi(r)=d\,{\rm arctg}\,\frac{\widetilde{y}(r)}{\widetilde{x}(r)}=\frac{1}{1+\Bigl(\frac{\widetilde{y}(r)}{\widetilde{x}(r)}\Bigr)^2}
\left(\frac{\widetilde{y}\,'(r)}{\widetilde{x}(r)}-\frac{\widetilde{y}(r)\widetilde{x}\,'(r)}{\big(\widetilde{x}(r)\big)^2} \right)dr\\
=\widetilde{x}\,'(r)\,\cos^2\phi(r)\,\frac{\frac{\widetilde{y}\,'(r)}{\widetilde{x}\,'(r)}-{\rm tg}\,\phi(r)}{\widetilde{x}(r)} \,dr=
\widetilde{x}\,'(r)\,\cos\phi(r)\,\frac{{\rm tg}\,(\vartheta_t-\vartheta_{\xi})-{\rm tg}\,\phi(r)}{r}\,dr\,.
\end{multline*}

Note that for a smooth arc $\widetilde{\gamma}$, for each $r\in(0,r_0]$ there exists $r_*\in[0,r]$ such that for
$t_*=\xi+r_*\,e^{i(\phi(r_*)+\phi_0)}$ the following relations are fulfilled:
\begin{equation}\label{spiv-arg-dotych}
 |\phi(r)|=|\vartheta_{t_*}-\vartheta_{\xi}|\le \omega_{\gamma}(\vartheta_t,r).
\end{equation}

Without loss of generality, we consider $r_0$ to be small enough to satisfy the inequality $\omega_{\gamma}(\vartheta_t,r_0)<1$.
Then we get the estimate
\[\int\limits_{\widetilde{\gamma}}|d\arg(t-\xi)|\le\, c\, \int\limits_{0}^{r_0}\frac{\omega_{\gamma}(\vartheta_t,r)}{r}\,dr<\infty, \]
where the constant\, $c$\, depends on $r_0$, but does not depend on $\xi$, generally speaking.

Finally, using the obtained estimate, we estimate the variation
\begin{equation*}
V_{\gamma}[\arg(t-\xi)]\le\int\limits_{\gamma[t_1,\xi]}|d\arg(t-\xi)|+\int\limits_{\gamma[\xi,t_2]}|d\arg(t-\xi)|+
\int\limits_{\gamma\setminus\gamma_{r_0}(\xi)}\frac{|dt|}{|t-\xi|}
\le 2c\, \int\limits_{0}^{r_0}\frac{\omega_{\gamma}(\vartheta_t,r)}{r}\,dr+\frac{{\rm mes}\,\gamma}{r_0}\,,
\end{equation*}
which yields the fulfillment of the condition \eqref{obm-var-arg} for the curve  $\gamma$.
\hfill $\Box$

\vskip 2mm

It is obvious that the condition \eqref{um-Dini-dotychna} is a weaker constraint on the curve $\gamma$ compared to the condition \eqref{um-Dini-ln-dotychna}.

The relation \eqref{spiv-arg-dotych} implies the estimate
\[ \omega_{\widetilde{\gamma}}\big(\arg (t-\xi),\eta\big)\le \omega_{\gamma}(\vartheta_t,\eta) \qquad \forall\,\eta\in [0,r_0],  \]
where the arc  $\widetilde{\gamma}$ is defined in the proof of Proposition \ref{prop-3}.
Therefore, if the modulus of continuity of the angle $\vartheta_t$ satisfies the Dini condition \eqref{um-Dini-dotychna}, then the condition of the same form is also satisfied for the modulus of continuity $\omega_{\widetilde{\gamma}}\big(\arg (t-\xi),\eta\big)$ at all points $\xi\in\gamma$:
\begin{equation}\label{um-Dini-arg}
 \int\limits_{0}^{1}\frac{\omega_{\widetilde{\gamma}}\big(\arg (t-\xi),\eta\big)}{\eta}\,d\eta<\infty.
\end{equation}

Let us show that the class of smooth Kr\'al curves differs from the class of smooth curves $\gamma$ that satisfy the conditions of the form \eqref{um-Dini-arg} at all points $\xi\in\gamma$. First, we give an example of a smooth curve $\gamma$ that is a Kr\'al curve, but the condition  \eqref{um-Dini-arg} is not satisfied at a point $\xi\in\gamma$.

\vskip 1mm

\begin{example}\label{ex-2}
Consider the smooth arc
\[\widetilde{\gamma}=\left\{t(r)=r\exp\Big(-i\,\frac{1}{\ln r}\Big) :  r\in(0,r_0]\right\}, \]
where $r_0$ is the smallest positive root of the equation ${\rm Re}\,t'(r)=0$.
It is obvious that the one-sided tangent to the arc $\widetilde{\gamma}$ at the beginning point
$t_0=0$ is the positive semi-axis of the real axis.
At the end point $t(r_0)$, the arc $\widetilde{\gamma}$ has the one-sided tangent parallel to the imaginary axis of the complex plane.

Let $\Gamma$ be such an arc of the ellipse that includes the points $z=x+iy$ satisfying the equation
\[\frac{x^2}{({\rm Re}\,t(r_0))^2}+\frac{(y-{\rm Im}\,t(r_0))^2}{({\rm Im}\,t(r_0))^2}=1, \]
which is smoothly glued to the arc $\widetilde{\gamma}$ at the points $0$ and $t(r_0)$.
Then $\gamma=\widetilde{\gamma}\cup\Gamma$ is a closed smooth Jordan curve.

It is obvious that the curve $\gamma$ satisfies the condition \eqref{obm-var-arg} because
$V_{\gamma}[\arg(t-\xi)]=\pi$ for all $\xi\in\gamma$.
At the same time, the condition \eqref{um-Dini-arg} is not satisfied at the point $\xi=0$ owing to the fact that
\[\int\limits_{0}^{1}\frac{\omega_{\widetilde{\gamma}}\big(\arg t,\eta\big)}{\eta}\,d\eta\ge
-\int\limits_{0}^{r_0}\frac{1}{\eta\ln \eta}\,d\eta=\infty.\]
\end{example}

\vskip 2mm

Now we give an example of a smooth curve $\gamma$ for which the conditions of the form \eqref{um-Dini-arg} are satisfied at all points $\xi\in\gamma$, but it is not a Kr\'al curve.

\vskip 2mm

\begin{example}\label{ex-3}
Consider the smooth arc
\[\Gamma_1=\left\{t(r)=r\exp\Big(-i\,\frac{r}{\ln r}\,\cos\frac{\pi}{r}\Big) :  r\in(0,1/2]\right\}. \]

Let $\Gamma_2$ be such an arc of the ellipse that includes the points $z=x+iy$ satisfying the equation
\[\frac{x^2}{a^2}+\frac{(y-b)^2}{b^2}=1 \]
with fully defined positive $a$ and $b$,
which is smoothly glued to the arc $\Gamma_1$ at the points $0$ and $t(1/2)$.
Then $\gamma=\Gamma_1\cup\Gamma_2$ is a closed smooth Jordan curve.

For each point  $\xi\in\gamma$, consider the arcs $\gamma[t_-,\xi]$ and $\gamma[\xi,t_+]$ defined in the proof of Proposition \ref{prop-3}, and denote them by $\gamma_{\xi}^-$ and $\gamma_{\xi}^+$, respectively. Considering the function $\arg (t-\xi)$ on the arc $\gamma_{\xi}^{\pm}$, we redefine it at the point $t=\xi$ by the limiting value
\[\lim\limits_{t\to\xi, t\in\gamma_{\xi}^{\pm}}\arg (t-\xi).\]

As a result, for each $\xi\in\gamma$, the function $\arg (t-\xi)$ satisfies the H\"older condition on each of the arcs
$\gamma_{\xi}^-$ and $\gamma_{\xi}^+$:
\[ |\arg (t_1-\xi)-\arg (t_2-\xi)|\le c\,|t_1-t_2|^{\alpha} \qquad \forall\,t_1,t_2\in \gamma_{\xi}^{\pm} \]
for all $\alpha\in(0,1/2]$, where the constant\, $c$\, does not depend on $t_1$ and $t_2$. Therefore,
the conditions of the form \eqref{um-Dini-arg} are satisfied at all points of  $\xi\in\gamma$.

At the same time,
\[V_{\gamma}[\arg t]\ge V_{\Gamma_1}[\arg t]\ge \sum\limits_{n=2}^{\infty} \frac{1}{n\ln n}=\infty,\]
i.e., the condition \eqref{obm-var-arg} is not satisfied for the curve $\gamma$.
\end{example}

\vskip 2mm

Example \ref{ex-3} also shows that the condition  \eqref{um-Dini-dotychna} on the angle $\vartheta_t$ in Proposition \ref{prop-3} can not be replaced by similar conditions of the form \eqref{um-Dini-arg} on the function $\arg(t-\xi)$.

\section{Sufficient conditions for the continuous extension of the real part of the Cauchy-type integral to the boundary of domain with unbounded variation of the function\, {\boldmath\mbox{$\arg(t-\xi)$}}}
\label{Sufficient conditions}

We shall now consider curves for which the condition \eqref{obm-var-arg} is not satisfied, generally speaking.

In what follows, we use the following characteristic of the function $f \colon E\rightarrow \mathbb{C}$ continuous on the set $E\subset \mathbb{C}$ (see O.F.~Gerus \cite{G96-2-4}):
\[ \Omega_{E,f}(a,b):=\sup\limits_{a\le \eta\le b}\frac{\omega_{E}(f,\eta)}{\eta} \qquad \mbox{for} \quad 0<a\le b. \]
The function $\Omega_{E,f}(a,b)$ does not increase monotonically with respect to the variable $a$ and does not decrease monotonically with respect to the variable $b$. In addition, the function $a\,\Omega_{E,f}(a,b)$ does not decrease monotonically with respect to the variable $a$.

Denote\, $E^{R,\psi_1,\psi_2}(\xi):=\{z=\xi+re^{i\phi} : R/2<r<R, \psi_1<\phi<\psi_2\}$.

Let us describe a certain finite set of Jordan arcs placed in the closure of domain $E^{R,0,\psi}(0)$.
For this purpose, we consider two sets of points $\{\tau_j\}_{j=1}^n$ and $\{\eta_j\}_{j=1}^n$ located on the rectilinear parts of the boundary of domain $E^{R,0,\psi}(0)$ such that
\[R\ge\tau_1\ge\tau_2\ge\dots\ge\tau_n\ge R/2, \qquad \eta_j=|\eta_j|\,e^{i\psi} \quad \mbox{and} \quad
R\ge|\eta_1|\ge|\eta_2|\ge\dots\ge|\eta_n|\ge R/2.\]

Let $\Gamma:=\bigcup_{j=1}^{n} \Gamma_j$, where $\Gamma_j$ is a Jordan arc with ends at the points  $\tau_j$ and $\eta_j$. Moreover, the arcs
$\Gamma_j$, $j=1,2,\dots,n$, excluding the ends, lie in the domain  $E^{R,0,\psi}(0)$ and
pairwise do not intersect within this domain.
In addition, if the arc $\Gamma_j$ is oriented from the point $\tau_j$ to the point  $\eta_j$, then the next arc $\Gamma_{j+1}$ is oriented in the opposite direction from the point  $\eta_{j+1}$ to the point $\tau_{j+1}$, and vice versa,
if the arc $\Gamma_j$ is oriented from the point $\eta_j$ to the point $\tau_j$, then the arc $\Gamma_{j+1}$ is oriented from the point $\tau_{j+1}$ to the point $\eta_{j+1}$.

Consider the auxiliary statements. The lemmas given below are proved by means of the method developed by T.S.~Salimov \cite{Salimov} for estimating the modulus of continuity of the Cauchy singular integral on an arbitrary closed rectifiable Jordan curve and adapted by O.F.~Gerus \cite{G96-2-4} for the purpose of using the modulus of continuity of the integral density instead of the regularized (by means of the Stechkin construction) modulus of continuity, which is used in the paper \cite{Salimov}.

\vskip 1mm

\begin{lemma}\label{lem-1}
{\it If a function  $f : \Gamma\rightarrow\mathbb{R}$ is continuous on $\Gamma$, then
\[\bigg|\,\int\limits_{\Gamma} f(t)\,d\arg t\,\bigg|\le\, \bigg( R\,\Omega_{\Gamma,f}\Big(\frac{R}{n},R\Big)+
\max\limits_{t\in\Gamma} |f(t)|\bigg)\,\psi+\frac{2\,\omega_{\Gamma}(f,\lambda)\,{\rm mes}\,\Gamma}{R}\,,\]
where $\lambda:=\max\limits_{j}\,{\rm mes}\,\Gamma_j$ and\, $\arg t$\, is any branch of the multivalued function
${\rm Arg}\,z$, which is continuous on~$\Gamma$.}
\end{lemma}

\vskip 1mm

For a given closed rectifiable Jordan curve $\gamma$, we shall consider its intersections with the domains $E^{R,\psi_1,\psi_2}(\xi)$, where $\xi\in\gamma$. We shall denote these intersections by  $\gamma_{R,\psi_1,\psi_2}(\xi)$.
By $n_{\gamma}(\xi,R,\psi_1,\psi_2)$ we denote the number of connected components of the set
$\gamma_{R,\psi_1,\psi_2}(\xi)$, the ends of which lie on the different segments $\{z=\xi+re^{i\psi_1} : R/2\le r\le R\}$ and $\{z=\xi+re^{i\psi_2} : R/2\le r\le R\}$.
We can say that the number $n_{\gamma}(\xi,R,\psi_1,\psi_2)$ expresses the number of complete oscillations of the function $\arg(t-\xi)$ in the domain $E^{R,\psi_1,\psi_2}(\xi)$.
Since the curve $\gamma$  is rectifiable, the number $n_{\gamma}(\xi,R,\psi_1,\psi_2)$ is finite, but with fixed $\xi$ and $R$ it can tend to infinity when $\psi_2\to\psi_1$.

Consider the case where there exist $\xi\in\gamma$ and $R\in(0,d]$ such that
\begin{equation}\label{pok-kolyv}
k_{\gamma}(\xi,R):=\max\, \Big\{1, \,\, \sup\limits_{0\le\psi_1<\psi_2<2\pi} n_{\gamma}(\xi,R,\psi_1,\psi_2)\Big\} <\infty.
\end{equation}

Denote by $\varphi_{\gamma}(\xi,R)$ the Lebesgue measure (given on the segment $[0,2\pi]$) of the set of those $\phi\in[0,2\pi]$ for which the rays $\{z=\xi+re^{i\phi} : r>0\}$ have a nonempty intersection with the set $\gamma_{R}(\xi)\setminus\gamma_{R/2}(\xi)$.

\vskip 2mm

\begin{lemma}\label{lem-2}
{\it Let a closed Jordan curve $\gamma$ be Ahlfors-regular and satisfy the condition \eqref{pok-kolyv} and let a function
$g \colon \gamma\rightarrow \mathbb{R}$ be continuous on $\gamma$.
Then the following estimation holds:}
\[\bigg|\,\int\limits_{\gamma_{R}(\xi)\setminus\gamma_{R/2}(\xi)} \big(g(t)-g(\xi)\big)\,d\arg (t-\xi)\,\bigg|
\le\, 6 \,
R\,\varphi_{\gamma}(\xi,R)\,\Omega_{\gamma,g}\Big(\frac{R}{k_{\gamma}(\xi,R)}\,,\,R\Big).\]
\end{lemma}

\vskip 2mm

\begin{lemma}\label{lem-3}
{\it Let a closed Jordan curve $\gamma$ be Ahlfors-regular and let a function
$g \colon \gamma\rightarrow \mathbb{R}$ be continuous on $\gamma$.
Let for $\xi\in\gamma$ the condition \eqref{pok-kolyv} be satisfied for all $R\in[\delta,2\varepsilon]$, where
$0<\delta<\varepsilon\le d/2$. Then the following estimation holds:
\begin{equation}\label{main-lem-oz}
\bigg|\,\int\limits_{\gamma_{\varepsilon}(\xi)\setminus\gamma_{\delta}(\xi)}
 \big(g(t)-g(\xi)\big)\,d\arg(t-\xi)\,\bigg|
 \le\, c\, \bigg(\int\limits_{\delta}^{2\varepsilon}\, \widehat{\varphi}_{\gamma}(\xi,\eta)\,\Omega\,_{\gamma,\,g}\Big(\frac{\eta}{\widehat{k}_{\gamma}(\xi,\eta)}\,,\,\eta\Big)\,d\eta+ \omega_{\gamma}(g,\varepsilon)\bigg),
 \end{equation}
where\,\,\, $\widehat{\varphi}_{\gamma}(\xi,R):=\sup\limits_{r\in[R/2,\,R]}\, \varphi_{\gamma}(\xi,r)$\,, \,\, $\widehat{k}_{\gamma}(\xi,R):=\sup\limits_{r\in[R/2,\,R]}\, k_{\gamma}(\xi,r)$, and
the constant\, $c$\, does not depend on $\xi$, $\delta$ and $\varepsilon$.}
\end{lemma}

\vskip 2mm

The following statement is true:

\vskip 1mm

\begin{theorem}\label{theor-dost-oz}
{\it Let a closed Jordan curve $\gamma$ be Ahlfors-regular and let a function
$g \colon \gamma\rightarrow \mathbb{R}$ be continuous on $\gamma$.
Let there be a partition $\gamma=\gamma^1\cup\gamma^2$, for which there exists
$R_0\in(0,d]$ such that:

(a) for all $\xi\in\gamma^1$ and all $R\in(0,R_0]$, the curve $\gamma$ satisfies the condition \eqref{pok-kolyv} and, in addition, the following condition is satisfied:
\begin{equation}\label{um-Dini-comb}
\sup\limits_{\xi\in\gamma^1}\, \int\limits_{0}^{R_0} \widehat{\varphi}_{\gamma}(\xi,\eta)\,\Omega\,_{\gamma,\,g}\Big(\frac{\eta}{\widehat{k}_{\gamma}(\xi,\eta)}\,,\,\eta\Big)\,d\eta<\infty\,;
\end{equation}

(b) for each  $\xi\in\gamma^2$ there exists  $r(\xi)\in(0,R_0]$ such that the curve $\gamma$ satisfies the condition \eqref{pok-kolyv} for all
$R\in(r(\xi),R_0]$, the function $\arg(t-\xi)$ has bounded variation on the set $\gamma_{r(\xi)}(\xi)\setminus\{\xi\}$ and, in addition, the following condition is satisfied:
\begin{equation}\label{obm-var-arg+Dini}
\sup\limits_{\xi\in\gamma^2}\, \bigg(\,V_{\gamma_{r(\xi)}(\xi)}[\arg(t-\xi)]+
\int\limits_{r(\xi)}^{R_0} \widehat{\varphi}_{\gamma}(\xi,\eta)\,\Omega\,_{\gamma,\,g}\Big(\frac{\eta}{\widehat{k}_{\gamma}(\xi,\eta)}\,,\,\eta\Big)\,d\eta\bigg)
<\infty\,.
\end{equation}
Then the function ${\rm Re}\,\widetilde{g}(z)$ is continuously extended to the boundary $\gamma$  from the domain $D^+$ and $D^-$, and the limiting values $\big({\rm Re}\,\widetilde{g}\big)^{\pm}(\xi)$ are represented by the formulas \eqref{double-pot+} and \eqref{double-pot-} for all $\xi\in\gamma$.}
\end{theorem}

\vskip 1mm

\textit{\textbf{Proof.}} Let us show that under the conditions of theorem, the condition \eqref{nidu-nepr} is also satisfied.

Let $\varepsilon\in(0,R_0/2]$ and $\delta\in(0,\varepsilon)$. Then the estimation \eqref{main-lem-oz} holds for all
$\xi\in\gamma^1$, and for each $\xi\in\gamma^2$, taking into account Lemma \ref{lem-3}, we obtain the estimate
 \[\bigg|\,\int\limits_{\gamma_{\varepsilon}(\xi)\setminus\gamma_{\delta}(\xi)}
 \big(g(t)-g(\xi)\big)\,d\arg(t-\xi)\,\bigg|\le \bigg|\,\int\limits_{\gamma_{r(\xi)}(\xi)\setminus\gamma_{\delta}(\xi)}
 \big(g(t)-g(\xi)\big)\,d\arg(t-\xi)\,\bigg|\,\]
 \begin{multline*}
 +\bigg|\,\int\limits_{\gamma_{\varepsilon}(\xi)\setminus\gamma_{r(\xi)}(\xi)}
 \big(g(t)-g(\xi)\big)\,d\arg(t-\xi)\,\bigg|\\
  \le \omega_{\gamma}(g,\varepsilon)\,V_{\gamma_{r(\xi)}(\xi)}[\arg(t-\xi)]+
 c\, \bigg(\int\limits_{r(\xi)}^{2\varepsilon}\, \widehat{\varphi}_{\gamma}(\xi,\eta)\,\Omega\,_{\gamma,\,g}\Big(\frac{\eta}{\widehat{k}_{\gamma}(\xi,\eta)}\,,\,\eta\Big)\,d\eta+ \omega_{\gamma}(g,\varepsilon)\bigg),
 \end{multline*}
 where the constant\, $c$\, does not depend on $\xi$, $\delta$ and $\varepsilon$.

Under the conditions \eqref{um-Dini-comb} and \eqref{obm-var-arg+Dini}, the given estimates yield the condition \eqref{nidu-nepr}.

Now, to complete the proof, it remains to apply Theorem \ref{theor-crit}.
\hfill $\Box$

\vskip 1mm

\begin{corolary}\label{cor-1}
{\it The statement of Theorem \ref{theor-dost-oz} remains valid when the conditions \eqref{um-Dini-comb}, \eqref{obm-var-arg+Dini} are replaced by the conditions
\begin{equation}\label{um-Dini-comb-1}
\sup\limits_{\xi\in\gamma^1}\, \int\limits_{0}^{R_0} \frac{\widehat{\varphi}_{\gamma}(\xi,\eta)\,\widehat{k}_{\gamma}(\xi,\eta)\,\omega_{\gamma}(g,\eta)}{\eta}\,d\eta<\infty\,,
\end{equation}
\begin{equation}\label{obm-var-arg+Dini-1}
\sup\limits_{\xi\in\gamma^2}\, \bigg(\,V_{\gamma_{r(\xi)}(\xi)}[\arg(t-\xi)]+
\int\limits_{r(\xi)}^{R_0} \frac{\widehat{\varphi}_{\gamma}(\xi,\eta)\,\widehat{k}_{\gamma}(\xi,\eta)\,\omega_{\gamma}(g,\eta)}{\eta}\,d\eta\bigg)
<\infty\,,
\end{equation}
respectively.}
\end{corolary}

\vskip 2mm

It is clear that Corollary \ref{cor-1} follows from Theorem \ref{theor-dost-oz} and the inequality
\[\Omega\,_{\gamma,\,g}\Big(\frac{\eta}{\widehat{k}_{\gamma}(\xi,\eta)}\,,\,\eta\Big)\le \frac{\widehat{k}_{\gamma}(\xi,\eta)\,\omega_{\gamma}(g,\eta)}{\eta} \qquad
\forall\,\eta\in(0,R_0]. \]

Note that in the papers of T.S.~Salimov \cite{Salimov}, E.M.~Dyn'kin \cite{Dynkin} and O.F.~Gerus \cite{G96-2-4},
 the singular Cauchy integral \eqref{zved-siC} on an arbitrary closed rectifiable Jordan curve is considered under certain conditions on the integral density, which are reduced to a condition of the form \eqref{2-4:Dini} in the case of a curve satisfying the condition \eqref{2-4:nerivnist'-teta}.
Under such conditions on the curve and the integral density, the Cauchy-type integral \eqref{C-type-int} is continuously extended to the boundary $\gamma$ from the domains $D^+$ and $D^-$, and the limiting values $\big({\rm Re}\,\widetilde{g}\big)^{\pm}(\xi)$
are expressed by the formulas \eqref{double-pot+} and \eqref{double-pot-} for all $\xi\in\gamma$.

At the same time, under additional assumptions of the type \eqref{pok-kolyv} about the curve $\gamma$, Theorem
\ref{theor-dost-oz} and Corollary \ref{cor-1} allow to construct examples (see the next example) of the curves $\gamma$ that do not satisfy the condition \eqref{obm-var-arg} and the functions
$g$ that do not satisfy the condition \eqref{2-4:Dini}
and a similar condition (see N.A.~Davydov \cite{Davydov}) used in Theorem 2 in the paper of O.F.~Gerus and M.~Shapiro \cite{Ger-Sh-2}, but the limiting values
$\big({\rm Re}\,\widetilde{g}\big)^{\pm}(\xi)$ of the logarithmic double layer potential \eqref{log-pot} exist at all points
$\xi\in\gamma$ and are expressed by the formulas \eqref{double-pot+} and \eqref{double-pot-}.

\vskip 1mm

\begin{example}\label{ex-4}
Consider the curve
\begin{multline*}
\gamma=\Big\{z=e^{i\phi} : \phi\in[0,\pi]\Big\}\cup[-1,0]\cup\bigcup_{n=1}^{\infty}[2^{-2n+1},2^{-2n+2}]\\
\cup\bigcup_{n=1}^{\infty} \Big\{z=2^{-n}e^{i\phi} : \phi\in[0,1/n]\Big\}\cup \bigcup_{n=1}^{\infty} \Big\{z=re^{-i\,\frac{\ln 2}{\ln r}} : r\in[2^{-2n},2^{-2n+1}]\Big\}
\end{multline*}
and the function 
\[ g(t)=  \left \{ \begin{array}{lll}
                -1/\big(\ln |t|-1\big) & \mbox{for} & t\in\gamma\setminus\{0\},\\[1mm]
                0  & \mbox{for} & t=0
                \end{array} \right.
 \]
that does not satisfy the Dini condition \eqref{2-4:Dini}
as well as a similar condition used in Theorem 2 in the paper \cite{Ger-Sh-2} because
\[\min\left\{\int\limits_0^1\frac{\omega_{\gamma}(g,\eta)}{\eta}\,d\eta,\,\,\,\int\limits_{\gamma}\frac{|g(t)-g(0)|}{|t-0|}\,|dt|\right\}
\ge\int\limits_{-1}^{0}\frac{|g(t)-g(0)|}{|t-0|}\,dt=-
\int\limits_{-1}^{0}\frac{dt}{|t|\,(\ln |t|-1)}=\infty\,.\]

For $0<\varepsilon\le 1/2$, denoting by $n_0$ the smallest natural number $n$ that satisfies the inequality $2^{-n}\le\varepsilon$, we obtain the estimate
\begin{equation}\label{gamma-Alf-cur}
 \theta(\varepsilon)=\theta_{0}(\varepsilon)\le 2\varepsilon+\sum\limits_{n=n_0}^{\infty}\frac{1}{n\,2^n}+\int\limits_{0}^{\varepsilon}\sqrt{1+\frac{\ln^2 2}{\ln^4 r}}\,dr \le 2\varepsilon+2\varepsilon+\frac{\sqrt{1+\ln^2 2}}{\ln 2}\,\varepsilon\le \left(4+\frac{\sqrt{1+\ln^2 2}}{\ln 2}\right)\varepsilon\,,
\end{equation}
which proves the fulfillment of the condition \eqref{2-4:nerivnist'-teta} for the curve $\gamma$. Thus, $\gamma$ is an Ahlfors-regular curve.

At the same time,
\[V_{\gamma}[\arg t]\ge \sum\limits_{n=1}^{\infty} \frac{1}{n}=\infty\]
and the condition \eqref{obm-var-arg} is not satisfied for the curve $\gamma$, i.e., $\gamma$ is not a Kr\'al curve.

Let us show that the curve $\gamma$ and the function $g$ satisfy the conditions of Corolary \ref{cor-1}.
There is the partition $\gamma=\gamma^1\cup\gamma^2$, where
$\gamma^1=\{0\}$ and $\gamma^2=\gamma\setminus\{0\}$. Let $R_0=1/2$ and $r(\xi)=\min\,\{2|\xi|,1/2\}$ for all $\xi\in\gamma^2$.
Then for all $\xi\in\gamma^2$, we have the estimation
\[V_{\gamma_{r(\xi)}(\xi)}[\arg(t-\xi)]\le V_{\gamma_{r(\xi)/4}(\xi)}[\arg(t-\xi)]+
\int\limits_{\gamma_{r(\xi)}(\xi)\setminus\gamma_{r(\xi)/4}(\xi)}\frac{|dt|}{|t-\xi|}\le\pi+\frac{\theta\big(r(\xi)\big)}{r(\xi)/4}\le
\pi+4c,\]
where\, $c=4+(\ln 2)^{-1}\,\sqrt{1+\ln^2 2}$\, as it is follows from the estimate \eqref{gamma-Alf-cur}.

It is obvious that $\widehat{k}_{\gamma}(0,\eta)\le 2$ and $\widehat{\varphi}_{\gamma}(0,\eta)<-1/\ln\eta$ for all $\eta\in(0,R_0]$.

In addition, $\widehat{k}_{\gamma}(\xi,\eta)\le 2$ for all $\xi\in\gamma^2$ and all $\eta\in(r(\xi),R_0]$.

Finally, taking into account the relations $|t|\le |t-\xi|+|\xi|\le 3\eta/2$, which hold for all $\xi\in\gamma^2$ and all $t\in\gamma$ such that
$|t-\xi|=\eta\in(r(\xi),R_0]$, for the specified $\xi$ and $\eta$, we obtain the inequality
$\widehat{\varphi}_{\gamma}(\xi,\eta)<-3/\ln(3\eta/2)$.

Now, the fulfillment of the conditions \eqref{um-Dini-comb-1} and \eqref{obm-var-arg+Dini-1} follows from the estimation
\[\sup\limits_{\xi\in\gamma}\,\int\limits_{\beta(\xi)}^{1/2} \frac{\widehat{\varphi}_{\gamma}(\xi,\eta)\,\widehat{k}_{\gamma}(\xi,\eta)\,\omega_{\gamma}(g,\eta)}{\eta}\,d\eta<
6\,\int\limits_{0}^{1/2} \frac{d\eta}{\eta\,\ln(3\eta/2)(\ln\eta-1)}<\infty\,, \]
where $\beta(\xi)=0$ for $\xi=0\in\gamma^1$ and $\beta(\xi)=r(\xi)$ for $\xi\in\gamma^2$.

Thus, all conditions of Corolary \ref{cor-1} are satisfied for the given curve $\gamma$ and the given function $g$.

As a result, we can state that the limiting values
$\big({\rm Re}\,\widetilde{g}\big)^{\pm}(\xi)$ of the logarithmic double layer potential \eqref{log-pot} exist at all points $\xi\in\gamma$ and are expressed by the formulas \eqref{double-pot+} and \eqref{double-pot-}.
\end{example}

\vskip 2mm

{\bf Acknowledgements.} The author acknowledges to Prof. Massimo Lanza de Cristoforis and Prof. Sergei Rogosin for the help with finding some references and Prof. Oleg Gerus for very useful discussions of results.
The author acknowledges also the financial support of UNIPD SRF and INdAM.

\end{document}